\documentclass[]{interact}

\usepackage{epstopdf}% To incorporate .eps illustrations using PDFLaTeX, etc.

\usepackage{natbib}% Citation support using natbib.sty
\bibpunct[, ]{(}{)}{;}{a}{}{,}% Citation support using natbib.sty
\renewcommand\bibfont{\fontsize{10}{12}\selectfont}% Bibliography support using natbib.sty

\theoremstyle{plain}% Theorem-like structures provided by amsthm.sty

\theoremstyle{definition}

\theoremstyle{remark}

\begin{document}

\articletype{ARTICLE TEMPLATE}

\title{Optimising structure in a networked Lanchester Model for Fires and Manoeuvre in Warfare}

\author{
\name{
Alexander C. Kalloniatis\textsuperscript{a}\thanks{CONTACT A.~C. Kalloniatis. Email: alexander.kalloniatis@dst.defence.gov.au},
Keeley Hoek\textsuperscript{b},
Mathew Zuparic\textsuperscript{a} and
Markus Brede\textsuperscript{c}
}
\affil{
\textsuperscript{a} 
Joint and Operations Analysis Division,
Defence Science and Technology Group, Canberra, ACT, Australia;
\textsuperscript{b}
Mathematical Sciences Institute, Australian National University,
Canberra, Australia;
\textsuperscript{c} 
Agents, Interactions and Complexity Group,
Department of Electronics and Computer Science, University of Southampton, Southampton, UK.
}
}

\maketitle
\begin{abstract}
We present a generalisation of the classical Lanchester model for directed
fire between two combat forces but now employing networks for the 
manoeuvre of Blue and Red forces, and the pattern of engagement between the two. The model therefore 
integrates fires between dispersed elements, as well as manoeuvre through an
internal-to-each-side diffusive interaction. We explain the model with several simple examples, including cases where conservation laws hold. We then apply an optimisation approach where, for a fixed-in-structure adversary, we optimise the internal manoeuvre and external engagement structures where the trade-off 
between maximising damage on the adversary and minimising own-losses can be examined. In the space of combat outcomes this leads to a sequence of transitions from defeat to stalemate and then to victory for the force with optimised networks. 
Depending on the trade-off between destruction
and self-preservation, the optimised networks 
develop a number of structures including
the appearance of so-called sacrificial nodes,
that may be interpreted as feints, manoeuvre hubs,
and suppressive fires. 
We discuss these in light of Manoeuvre Warfare theory.
\end{abstract}

\begin{keywords}
Warfighting; mathematical model; dynamics; networks; manoeuvre; optimisation
\end{keywords}

\section{Introduction}
In the age of high powered computers it has become accepted folklore
that the days of understanding warfare through equation-based models are over. This is
argued to be the case because of the multiplicity of factors influencing modern warfare, such as dispersed forces, networked Command and Control (C2), diverse rivals and stakeholders in conflicts, and technologies for Intelligence-Surveillance-Reconnaissance (ISR). In this paper we propose
and explore a model that offers a platform to bring such factors
together while still retaining the advantages of equation based approaches.
At the heart of this remain the Lanchester equations for combat\cite{Lan1916}. For directed fires these are
\begin{eqnarray}
\dot{B}(t) &=& \gamma_B B(t) - \kappa_R R(t) \nonumber \\
\dot{R}(t) &=& \gamma_R R(t) - \kappa_B B(t)
\label{LanchDirF}
\end{eqnarray}
where $B(t),R(t)$ represent the total number of forces of Blue and Red
at time $t$, $\gamma_B, \gamma_R$ the rate of resupply, and $\kappa_B,\kappa_R$ represent their respective kill-rates of the other. The $\kappa$ are also variously referred to
as fire-rates or fire-power or lethality, terms which we
shall use interchangeably.
The so-called `undirected fire' version of these have $B(t)R(t)$
on the right-hand side of the equations.
Proposed in 1916 by F.W. Lanchester, these simple equations have
seen diverse uses beyond the original context of the 
statically arrayed massed fires
of World War I. Such extensions will be reviewed below, but
nascent is the incorporation of 
recent concepts such as as Network Centric Warfare \cite{Alberts1999}, and realised in
General Stanley McChrystal's strategy for battling the Taliban in Afghanistan
``It takes a network to defeat a network"
\cite{McChrystal2011}. 
Nevertheless, a fully networked version of
the Lanchester model has not, to our knowledge, been written down.

In other dynamical systems, the network approach is  well-developed, offering insights into various complex systems -- biological, social and technological
\cite{BornSchus2003} -- through numerical simulation and analytical approximation.
The network paradigm has advanced to {\it multi-layer} networks
\cite{BoccEtAl2014}, where
the connectivity within any layer
and across layers may vary considerably.
In warfare such layers may describe how forces are internally distributed, how they engage with the enemy, how C2 is
arranged, and how ISR
systems enable situation awareness. 
This paper takes initial
steps towards capturing some of these using
the multi-layer formalism.

We propose a form of
the Lanchester equations whereby
rival heterogeneous forces manoeuvre
internal resources, and
apply force against each other on external networks. We explore optimisation of these networks and
observe emergent behaviours that may
be interpreted through known
manoeuvre warfare concepts.

We first review the literature on the Lanchester model,
situating our formulation in recent
variations of the model.
We propose our networked Lanchester model
and examine a small scale example to illustrate its dynamics.
We then outline
the optimisation and explore the trade-off between
a force seeking to minimise its losses versus
maximising its degradation of the adversary.
In particular, we observe
behaviours that may be interpreted
as `concentration', `suppressive fire' and `feints', concepts
well-known in manoeuvre warfare. 
We show how some of these behaviours arise analytically exploiting
the power of the Lanchester analytical framework. We then conclude,
discuss these aspects in light of the literature on
Manoeuvre Theory of warfare, 
and offer prospects for the future. An appendix provides
analytical explanation for some of the numerical results
and supplementary material presents an alternate constraint in the
optimisation scenario.

\section{Status of Lanchester combat modelling}
Enhancements of the original Lanchester equations are many.
As a description of the waxing and waning of two abstract
entities, the equations may be interpreted more generally than
mortal combat, for example in terms of network
cyber-attack and defence \cite{LiuEtAl2013}.
Mathematical generalisations include: general power law $B^\alpha R^\beta$ forms for the right hand side \cite{Epstein1997} and other nonlinear forms
\cite{KimEtAl2017}, stochastic diffusive models with advection \cite{Protopopescu1989} and swarming effects \cite{Keane2011} as partial differential equations, and game theory
applications \cite{Hohzaki2016}. These still only represent two homogeneous sets
of forces arrayed against each other. In
this era of joint or combined operations and multi-role platforms,
the need for representing mixed-forces becomes acute.
For such mixed-forces, numerous authors have proposed the $(N,M)$ Lanchester problem of $N$ Blue forces
arrayed against $M$ Red forces. Suppressing, at this stage,
logistics/replenishment, this can be expressed through the
system ($i=1, \dots, N; j=1,\dots, M)$
\begin{eqnarray}
\dot{B}_i(t) = - \sum_{j=1}^M \mu(B)_{ij} R_j \nonumber \\
\dot{R}_j(t) = - \sum_{i=1}^N \mu(R)_{ji} B_j, 
\label{mixedLan}
\end{eqnarray}
where the $\mu$ capture the distribution of fire between units of Blue and Red respectively
in a manner that tailors the rate to the adversary unit; for this reason
some functional dependence of the $\mu$ on $B$ or $R$ is indicated. (Note that
many authors use upper and lower indices according to attack and defence,
which we will treat symmetrically).

Early such approaches are\cite{Rob1992} on the $(2,1)$ problem and
\cite{Cole1993} on the $(2,2)$ case who typically
used time-independent $\mu(B)_{ij},\mu(R)_{ij}$. 
It was pointed out in\cite{KaupEtAl2005}
that such heterogeneous models cannot be 
interpreted in terms of the reality of combat;
significant discontinuities in the redistribution of fire must be taken into account when units reach the value zero,
when combat entities die. They 
proposed a model where the lethality depends
on the attacker and the defender such that
when one defender is destroyed the attacker may
engage other targets.

These concerns have been furthered in a series
of works authored solely or jointly by MacKay, starting with
\cite{MacKay2009} where analytic solubility could be gained
by considering certain separable 
but dynamical forms $\mu_{ij}(t)=\kappa_i \mu_j(t)$,
using weighted-ratios of one force element to the total
weighted-sum of the
force, for example 
$\mu_i(t)=\rho_i B_i(t)/\sum_j \rho_j B_j(t)$.
These reflect that a commander may have
a choice in how to distribute fire at the start 
of a battle but thereafter it evolves according to
the proportion of units of
a particular type. This was
challenged by \cite{LiuEtAl2012}
who allowed for targets to change during
battle as adversary units were eliminated. This also was
questioned by MacKay, 
noting that in the fog-and-friction of war,
such reallocation of targets may be unrealistic
after the battle has begun\cite{MacKay2012}.
We return to these points when we develop our model.

MacKay's formulations generate conserved quantities generalising the
Lanchester square law for directed fire, $\kappa_R R(t)^2 - \kappa_B B(t)^2$ 
when $\gamma_B=\gamma_R=0$ in Eqs.(\ref{LanchDirF}).
This determines the victor based on initial conditions. All
these authors recognise very well
that these equations purport to describe
battle of annihilation; quantities changing sign have no 
meaning. Within a `kinetic effects' interpretation
(as opposed to non-lethal, such as
cyber \cite{LiuEtAl2013}) the mixed forces models remain representations of unmitigated attrition, of ``perfectly horrible" warfare\cite{MacKay2009}.

\section{The networked Lanchester model}
Readers will recognise a network structure in Eqs.(\ref{mixedLan}):
the $\mu_{ij}$ may be interpreted as weighted adjacency matrix elements describing 
the engagement interaction between the heterogeneous
Blue and Red units. For $\mu_{ij}\in (1,0)$ we have a network of undirected links of nodes $(i,j)$ that are connected, respectively unconnected.
Though many authors
seek to model `networked forces' in Lanchester combat
\cite{TangLi2012},\cite{LiuEtAl2013}, none to our knowledge explicitly represent the
network structure in the dynamics of the combat.

This offers an opportunity of overcoming one of the limitations of
even the mixed Lanchester model, that through the course of battle
a force may be dynamically reconfigured to reinforce weak units or exploit weaknesses
of the adversary. This seeking to gain ``advantageous position relative to the enemy''
\cite{Lind1985} is known as `manoeuvre warfare', and is a further `warfighting function' after
attrition or `fires'.
`Manoeuvre' captures the idea that
success in the battlefield for forces with equally
matched firepower unit for unit may be gained by
not only providing more units,
but through skillful manipulation of those units.

We propose for this the diffusive interaction for Blue force units 
$i,j=1,\dots,N$:
\begin{equation}
\dot{B}_i(t) = \sum_{j=1}^N {\cal B}_{ij} \left(\delta_i^{(B)}(t)B_i(t) - \delta_j^{(B)}(t)B_j(t) \right)
\end{equation}
where $B_i(t) \in {\rm I\!R}$, ${\cal B}_{ij}$ is a $(0,1)$ adjacency matrix and $\delta_i^{(B)}(t)$ represent possibly
time-dependent weight factors between
different Blue force elements.
Such an approach has been used successfully for human
population migration \cite{Roman2017}.
It is trivially seen that for symmetric ${\cal B}_{ij}$
\begin{eqnarray}
\sum_{i=1}^N \dot{B}_i(t) = 0
\end{eqnarray}
as a consequence of the double sum becoming a difference 
under the sum of 
terms $k_i \delta_i B_i$, where we
use the convention of $k_i=\sum_j {\cal B}_{ij}=\sum_j{\cal B}_{ji}$
for the degree of node $i$. 

For the Red force, we introduce $R_i(t) \in {\rm I\!R}$, and
the corresponding internal adjacency matrix ${\cal R}_{lm}$
with nodes $l,m=1,\dots,M$ for Red agents. 
The networks ${\cal B}_{ij}$ and ${\cal R}_{lm}$
represent how a given distribution of resources $B_i(t) \geq 0$ and $R_l(t) \geq 0$ may {\it manoeuvre} dynamically through the course of battle;
we distinguish this from a future extension with {\it logistics} networks
that represent replenishment of the forces from outside the battle-space.
Thus ${\cal B}_{ij}$ and ${\cal R}_{lm}$ are denoted as
{\it manoeuvre} networks.
Characteristic constants $\gamma_B$ and $\gamma_R$
determine the strength or time-scale of transfer of forces within
the respective networks; these constants may be made node-dependent, not treated in this work.

Extending the network idea to the application of fires,
we have {\it engagement} networks 
${\cal E}^{(BR)}_{il}$ and ${\cal E}^{(RB)}_{li}$
(called $\mu_{ij}$ previously)
representing
the pattern of directed fires from one side to the other; these need not
be symmetrical between Blue and Red. Correspondingly kill-rates
$\kappa_B$ and $\kappa_R$, as with Eq.(\ref{LanchDirF}),
characterise the combat effectiveness of the respective forces.
These, like $\gamma_B, \gamma_R$, may also be node or link dependent.
Given the
complexity already of this model, we consider homogeneous forces
with uniform fire-power to understand the key 
impact of dynamical manoeuvre
within each side's forces. In the same spirit, we will use simplistic static engagement networks
in contrast to\cite{MacKay2009}.
Arguably, the ability to dynamically reallocate 
targets or the degree of fire-power may only be possible with
improved ISR and a functioning
C2 system able to integrate situation
awareness into decision-making. Indeed, as we shall argue further below,
sometimes in warfare strategies, are based on
the adversary's over-estimation of the strength of
some formations, or even that ``dead'' or ``dummy'' units
may indeed be the mis-directed focus of
attack.

Bringing these elements together gives an initial network generalisation of Eqs.(\ref{LanchDirF}):
\begin{eqnarray}
\dot{B}_i &=& -\gamma_B \sum_{j} {\cal B}_{ij} (\delta_i^{(B)} B_i - \delta_j^{(B)} B_j) 
-\kappa_R \sum_m {\cal E}^{(RB)}_{im} R_m \nonumber \\
\dot{R}_l &=& -\gamma_R \sum_{m} {\cal R}_{lm} (\delta_l^{(R)} R_l - \delta_m^{(R)} R_m)
-\kappa_B \sum_j {\cal E}^{(BR)}_{lj} B_j.
\label{initeqs}
\end{eqnarray}
This will not be the model's final form. The first terms in each
set of equations are the manoeuvre contributions, where resource
may be shifted through the respective networks according to relative
strengths. The second set of terms are the engagement contributions
but with predefined static distributions of fire.

Thus far, Eq.(\ref{initeqs}) describe a dynamics where the 
$B_i, R_l$ may become negative: the equations need factors
that forcing nodes to `drop out' once resource levels reach zero.
We may do this by Heaviside step functions $\Theta(x)$ or by a smeared
form using the hyperbolic tangent function
\begin{equation}
\Theta_\epsilon(x) \equiv \frac{1}{2} (1+ \tanh(x/\epsilon))
\end{equation}
which for small $\epsilon$ approximates a step function at $x=0$.
In the same spirit, the diffusive transfer of resources to a weaker
node should cease when that node has reached zero or a sufficiently
low level $\vartheta$. This may be achieved by using a shifted Heaviside
function, $\Theta(x-\vartheta)$ but for simplicity we set $\vartheta=0$.
The third modification of the equations is to moderate the attrition
term: in its form in Eq.(\ref{initeqs}) the entire resource of a node
may be brought to bear on multiple nodes of the adversary without
diminution.
Thus it is more realistic to divide the attrition by the out-degree
according to the engagement network,
\begin{equation}
k_i^{(RB)}\equiv \sum_m {\cal E}^{(RB)}_{im}
\end{equation}
for the engagement of Red with Blue agent $i$, with the same for $BR$.

Thus the final form of the equations reads:
\begin{eqnarray}
\dot{B}_i &=& - \gamma_B \sum_{j} {\cal B}_{ij} (\delta_i^{(B)} B_i - \delta_j^{(B)} B_j) \Theta_\epsilon(B_i)\Theta_\epsilon(B_j) \nonumber \\
&& -\kappa_R \sum_m {\cal E}^{(RB)}_{im}  \frac{R_m}{k_i^{(RB)}} \Theta_\epsilon(B_i) \Theta_\epsilon(R_m)\nonumber \\
\dot{R}_l &=& -\gamma_R \sum_{m} {\cal R}_{lm} (\delta_l^{(R)} R_l - \delta_m^{(R)} R_m) \Theta_\epsilon(R_l) \Theta_\epsilon(R_m)
\nonumber \\
&& -\kappa_B \sum_j {\cal E}^{(BR)}_{lj} \frac{B_j}{k_l^{(BR)}} \Theta_\epsilon(R_l)\Theta_\epsilon(B_j). \nonumber \\
\label{finaleqs}
\end{eqnarray}
For the weights $\delta$ we initially trialed constant values. Optimisation of the manoeuvre networks here generated
{\it disconnected} graphs with isolated strong nodes and poorly connected weak nodes of commensurate initial conditions. The explanation for
this is straightforward: connectivity here drains resource from strong nodes to weak which eventually die. 
Intuitively, channelling resource from strong to weak nodes
should take into account their relative strengths to
engaged adversaries.
We therefore propose the following form,
here for Blue's engagement with Red:
\begin{equation}
    \delta_i^{(B)} = \frac{1}{{\sum_m {\cal E}^{(BR)}_{im} R_m}+\epsilon},
    \label{delta-final}
\end{equation}
with some regularisation parameter $\epsilon$.

The choice of $\delta_i$ is a {\it warfighting heuristic}: not intrinsic to the model, it implements
a principle for how a force seeks to achieve manoeuvre.
While regularising divergences, $\epsilon$ also has physical meaning, noting that when
$ {\sum_m {\cal E}^{(BR)}_{im} R_m} $ vanishes,
$\epsilon$ scales into the rate $\gamma_R$. 
Thus, $\epsilon$ represents 
a form of `standing force' that a non-combatant node
seeks to maintain: for small $\epsilon$ non-engaged
nodes will give up resource rapidly to engaged
partners; for large $\epsilon$ all nodes in
the manoeuvre network will redistribute resource
regardless of the state of engagements.
In this work we will treat $\epsilon$ as infinitesimal.
The model of Eqs.(\ref{finaleqs}) offers a framework for
testing the consequences of such heuristics. Our purpose then  
is to show that such a manoeuvre heuristic
realises recognisable warfighting concepts,
so the validity of the model will be determined {\it post hoc}.

We will be interested in equally resourced,
fire-power matched but differently structured forces, 
so-called ``near peer adversaries'', as has become
the focus of recent shifts in the US \cite{FM3-0}
and UK \cite{Chuter2019} militaries, but without necessarily the
`massed force' paradigm \cite{Cole2019}. Thus, $B_i=R_i$ for engaged
units (the distinction between indices $i,j$ and $l,m$ now
becomes superfluous), or 
$\sum_j B_j=\sum_m R_m=T$.
Given this, and because the model is linear we may
divide through by the initial conditions $R_i(0)=B_i(0)$ or total force $T$. But from herein we
consider $0\leq B_i, R_l \leq 1$.

We will also set $\gamma_R=\gamma_B=1$
and vary the lethality $\kappa_R, \kappa_B$.
We will thus be examining
the impact of differing lethalities of
the two sides given different structures
for manoeuvre and engagement {\it for fixed time-scale
at which they are able to dynamically adapt} the
allocation of their force.

\section{Simple case study}
To illustrate the key dynamics of the model we consider
an engagement of two Blue units engaged against four Red units,
where the latter use a manoeuvre network to draw reserves
as defined in Eqs.(\ref{finaleqs},\ref{delta-final}).
\begin{figure}[tbp]
 \begin{center}
 \includegraphics[width=1.0\textwidth]{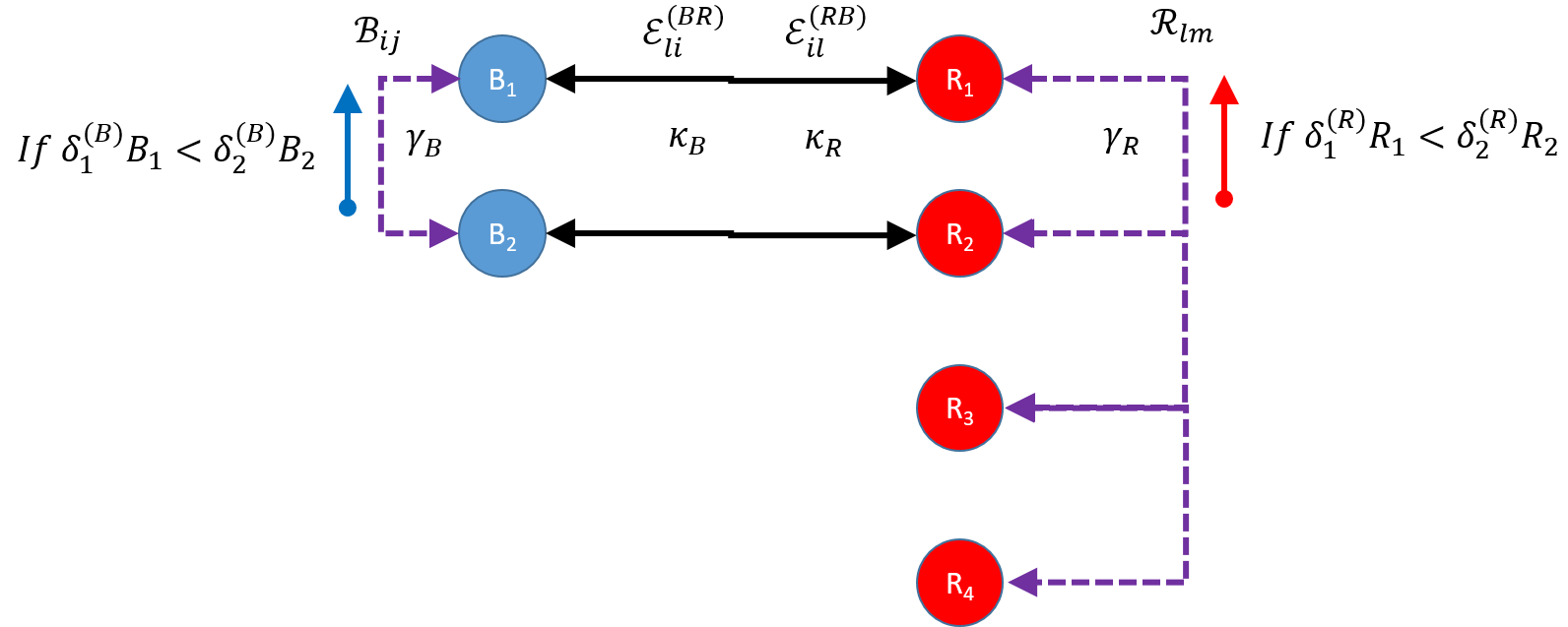}
  \caption {Illustration of the full model Eq.(\ref{finaleqs}) with two Blue against four Red, where solid black lines indicated
 directed fire and purple dashed lines indicate resource sharing paths. Network quantities are labelled across the top and rates below
 the top row. The directions of flow, indicated by blue and red
 arrows, are governed by the conditions involving the $\delta$;
 if the direction of the inequalities reverses so too does the direction
 flow.}
\label{casestudyB2R4}
\end{center}
 \end{figure}
We will consider a number of ways of determining where
manoeuvre provides advantage:
\begin{enumerate}
    \item initially equally matched adversaries, $B_i(0)=R_i(0)=1/2$ for
    $i=1,2$ but with Red having
    varied reserves in $R_3(0)=R_4(0)=f_R/2$, so that
    $\sum_j B_j(0)=1, \sum_m R_m(0)=1+f_R$;
    \item initially mismatched adversaries, $B_i(0)=1/2, R_i(0)=f_R/2 $ for
    $i=1,2$, but equal total, 
    $\sum_j B_j(0)=\sum_m R_m(0)=1$ so that $R_3(0)=R_4(0)=(1-f_R)/2$
    or zero if $f_R>1$; and
    \item initially mismatched adversaries, but where Red has extra
    reserves it can call upon, namely $B_i(0)=0.5, R_i(0)=f_R/2 $ for
    $i=1,2$, as well as $R_l(0)=f_R/2 $ for
    $l=3,4$ so that $\sum_j B_j(0)=1, \sum_m R_m(0)=2f_R$.
\end{enumerate}
With the Blue kill-rate fixed at $\kappa_B=1$, and equal manoeuvre
rates $\gamma_R=\gamma_B$ we vary both the Red kill-rate $\kappa_R$
and the fraction of the force available to Red.
We focus first on case $(3)$ for $\kappa_R<\kappa_B$
where Red has the most extreme initial disadvantage. How can Red's ability to manoeuvre spare resources
provide advantage?

\begin{figure}[tbp]
 \begin{center}
 \includegraphics[width=1.03\textwidth]{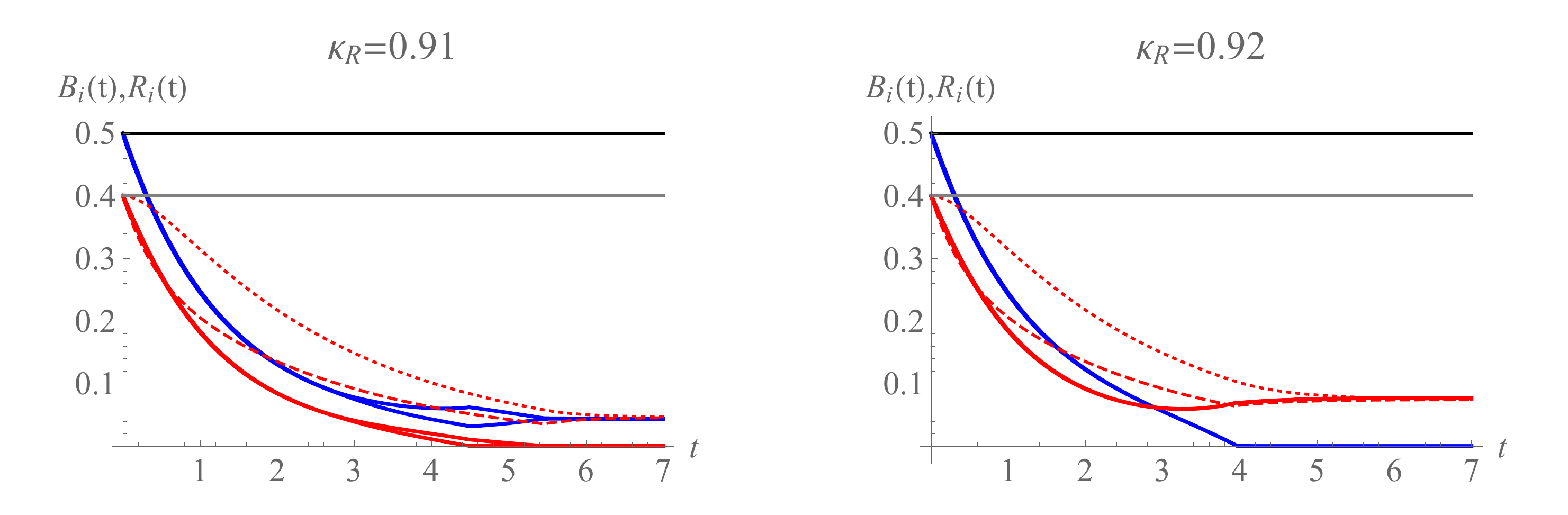}
  \caption {Case of two Blue against four Red dynamics for two different
  Red kill-rates, where Blue and Red solid lines represent engaged units and Blue and Red dashed/dotted lines represent supporting units. Black and gray represent the mean force per unit of Blue and Red respectively.}
\label{casestudyB2R4dyn}
\end{center}
 \end{figure}
 
In Fig.\ref{casestudyB2R4dyn} are the profiles for the 
units as functions of time for the case of $f_R=0.8$ and two
marginally different values of $\kappa_R=0.91, 0.92$ both less
than $\kappa_B=1$. Recall for case $(3)$ the total resources
available initially for Red are $1.6$ but each Red combat unit
begins with only $0.4$ compared to Blue's $0.5$.
We also show the mean initial resource per unit, black for Blue and
grey for Red so that Red is on average less resourced than
Blue.

For the lower value of $\kappa_R$ (left-hand panel), Blue's combat units maintain their
initial advantage; the dashed and dotted red lines
indicate Red reserves are offering resource to their combat partners,
but too slowly, given the combat units' insufficient fire-power, to save the day. When the Red combat units
have expired, their reserve counterparts stabilise as the combat has
finished. Blue has won in terms of the combat dynamics; the remaining
Red units are unable to transfer any more resource into the fight.

With a small increase of fire rate by Red, this situation is reversed. In the right hand plot of Fig.\ref{casestudyB2R4dyn} 
Red is able to manoeuvre 
resources into the fight with Blue, and Red combat units are able to
degrade Blue sufficiently, to enable Red to defeat Blue. After
Blue's expiry, Red reserve units continue to share resource with
combat units until they equalise and stabilise.

\begin{figure}[tbp]
 \begin{center}
 \includegraphics[width=0.8\textwidth]{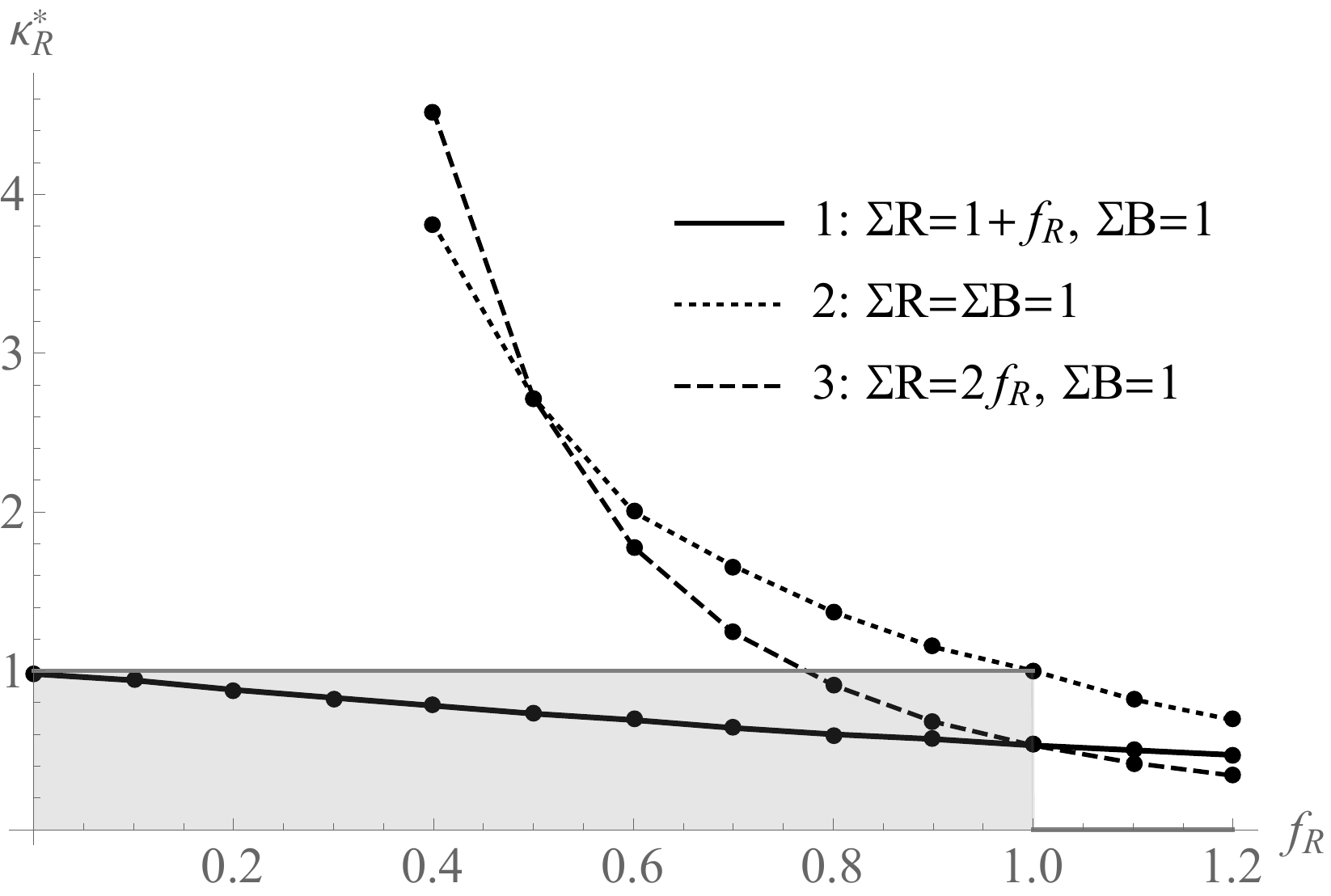}
  \caption {The critical value of $\kappa_R$ at which for
  a value of $f_R$ Red defeats Blue for the three different cases. With the grey region we indicate where
  Red gains advantage with lesser fire power and initial resource for its combat units.}
\label{casestudyB2R4critpt}
\end{center}
 \end{figure}

There is thus a balance between
reserve capacity, fire power and manoeuvering rate
enabling a set of combat units to defeat
a unit-for-unit superior force.

Within the same model we explore the
critical value of Red fire power, denoted $\kappa_R^*$,
at which Red is able to defeat Blue for the different cases
(1)-(3). To this end we numerically solve for sufficient
time until one or the other combat force is degraded
to zero for different $\kappa_R, f_R$ and for the three
cases of reserve supply. This is shown in Fig.\ref{casestudyB2R4critpt}. In the plot we
overlay a grey region where Red has weaker fire power
and initial combat unit strength compared to Blue.
In other words, curves that cross into the grey region
show where Red genuinely gains advantage significantly
because of its manoeuvre capability rather than
numerical or weapon strength.

For the second case with equal total resource
for both forces at the outset (dotted curve), $f_R$ must exceed the value one for Red to achieve advantage with lesser fire-power.
In other words, each Red combatant unit requires more
than half the initial resource of the Blue combatants.
Thus with the same total initial resource and lesser fire-power
and initial strength, the Red force is able to defeat Blue
through the manoeuvre of resource through its network.
However, with even more total initial resource,
it is possible for Red to achieve victory with
even lower fire-power and initial resource for
its combat units, as seen in the solid and dashed
curves for values inside the grey region.
We see then that the manoeuvre capability provides
advantage over and above pure
fire-power and strength.

\section{Optimal structures for manoeuvre and engagement}
We now amplify the scale of forces represented in
the model, and focus on how network heterogeneity
develops to provide advantage to one side even with less raw fire-power than the adversary. 
We will be interested
in the heterogeneity
that emerges according to the balance
of priorities of an optimal force.

\subsection{Network optimisation}
\label{networkopt}
We consider optimisation of battle outcomes for one side, here chosen as Red. 
We presume that Red's decision makers are interested in overall outcomes in which certain trade-offs $U_R(\lambda)$ between their own remaining forces 
$R(t)=\sum_i R_i(t)/N$ and the number of destroyed and remaining opponents $B(t)=\sum_i B_i(t)/N$ are optimised; here $N=\sum_i R_i(0)=\sum_i B_i(0)$
since all $B_i(0),R_i(0)=1, i=1,\dots,N$. 
Because of equal initial forces we cease distinguishing between
indices $(i,j,k,l) \in (1,\dots,N)$. We quantify utility for Red by
\begin{equation}
 U_R(\lambda) = \lambda R + (1-\lambda) (B(0)-B),
\end{equation}
where $B(0)$ gives the initial average force of Blue and $\lambda$ measures a trade-off between an interest in preserving their own forces
(`defence') and destroying the adversaries' forces (`offence'). 
Since we do not include resupply terms here this quantity will be non-negative.

The optimisation scheme is based on stochastic hill-climbing, to maximise $U_R$. In more detail: 
\begin{enumerate}
    \item Start with some initial engagement and manoeuvre networks $R_{ij},B_{ij}$ and ${\cal E}^{RB}_{ij}={\cal E}^{BR}_{ji}$,
based on Erdos-Renyi random graphs, 
which are built by randomly assigning links
between nodes according to a uniform probability;
here we use average degree
of the network $k=4$. The
combination of manoeuvre and engagement networks for one side we refer to as 
a `configuration'.
    \item Perform a random rewiring of Red's configurations by, either, with some probability $p$ rewiring Red's manoeuvre network $R_{ij}$ or, with probability $1-p$, rewiring of Red's engagement network 
    ${\cal E}^{RB}_{ij}$. For these modifications of the manoeuvre network, we pick a randomly selected link and move it to a randomly selected link vacancy, thus preserving the average degree of the manoeuvre network. 
    Similarly, for modifications in the engagement network, we select one engagement link uniformly at random and move it to a randomly selected  link vacancy in the engagement network. 
    Additionally, we also consider additions of new engagement links to randomly selected engagement network link vacancies or removal of randomly selected engagement links, such that Red can optimise the number of engagements it wants to participate in.
    \item For this modified structure we then numerically integrate the respective network Lanchester system Eqs.(\ref{finaleqs}) to determine stationary force concentrations and evaluate Red's utility $U_R(\lambda)$. New configurations are accepted if the stochastic modification resulted in an improvement in Red's utility, otherwise we reject the modification. In the first case, we keep the modified configuration and
    repeat (2). In the second case, we restore the previous battle structure before repeating step (2).
\end{enumerate}
To integrate the equations we use fourth order Runge-Kutta with time-step $\delta t=0.01$. Optimisation steps are repeated in the order of $10^5$ times, ensuring near-convergence to a final configuration. We repeat the procedure for a fixed number of different random initial conditions, and check final configurations to ensure the robustness of our findings presented below. We have also attempted simulated annealing to avoid trapping in local minima, however such localised
methods did not result in significant improvements compared to the method outlined here. One optimisation run of these $10^5$ iterations takes approximately 2 hours, but most of our results below have been averaged over 20 independent runs, resulting in a simulation time of approximately 2 days for one choice of parameters.

%In the inset to this
%panel we plot the density of surviving forces 
%of both sides at the stalemate point. Specifically,
%we observe that for small $L_R<1500$ we find a region (between the
%curves)  where both forces can be large and coexist, and for $L_A>1500$
%coexistence is only possible in a small region in which neither force is
%larger than $0.05$.

\begin{figure}[tbp]
\begin{center}
\includegraphics[width=.99\textwidth]{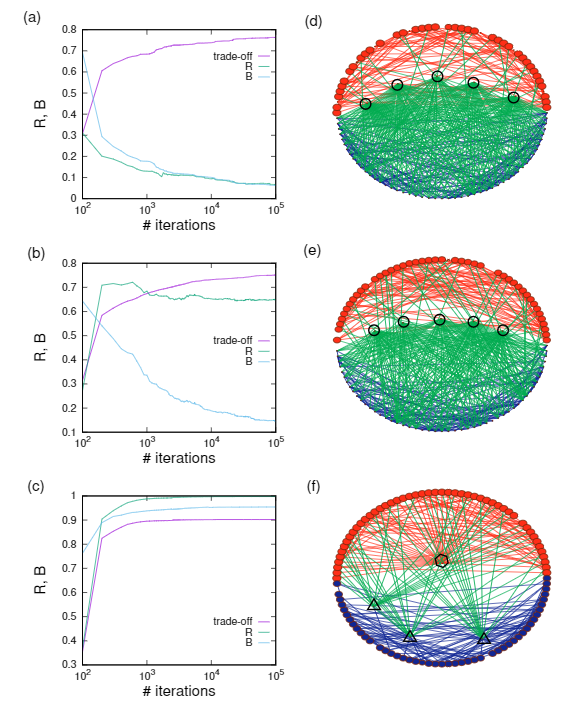}
 %\caption {(top row) Evolution of stationary values of forces and trade-off utility function for (a) low ($\lambda=0.2$), (b) intermediate ($\lambda=0.5$), and (c) large trade-off parameter $\lambda$ ($\lambda=0.9$). Trajectories are averaged over 20 independent runs. (middle row) Example networks at equilibrium after $10^5$ optimisation iterations for (d) low, (e) intermediate, and (f) large $\lambda$ as above. Red force nodes drawn in red, Blue force nodes in blue. Connections in the manoeuvre networks of Red and Blue forces in red or blue. Engagements are drawn in green.  (bottom row) Heat map of differences between normalised stationary Red and Blue forces vs. lethalities $\kappa_R$ and $\kappa_B$ for the best configuration found for networks optimized for $\kappa_R=0.5, \kappa_B=1$ and (g) $\lambda=0.5$ and (h) $\lambda=0.9$. Panel (i) shows the heatmap for a random  Blue and Red network with the numbers of links in the manoeuvre and engagement networks consistent with the optimised networks of Red. Red colours indicate regimes of larger Red surviving forces and Blue colours regimes of larger Blue surviving forces. In all cases the system size is $N_R=N_B=50$, $L_R=L_B=100$, $\kappa_R=0.5, \kappa_B=1$, and initial conditions $R_i=B_i=1$. }
 \caption {Evolution of stationary values of forces and trade-off utility function optimization for (a,d) low ($\lambda=0.2$), (b,e) intermediate ($\lambda=0.5$), and (c,f) large trade-off parameter $\lambda$ ($\lambda=0.9$). Trajectories are averaged over 20 independent runs. On the right, example networks 
 at equilibrium after $10^5$ optimisation iterations for (d) low, (e) intermediate, and (f) large $\lambda$ as above. Red force nodes drawn in red, Blue force nodes in blue. Connections in the manoeuvre networks of Red and Blue forces in red or blue. Engagements are drawn in green. Nodes surrounded by a black circle indicate red sacrificial nodes, blue nodes in a triangle denote blue nodes at which Red's attack is focused. The red node in a diamond in (f) indicates a manoeuvre hub. In all cases the system size is $N_R=N_B=50$, $L_R=L_B=100$, $\kappa_R=0.5, \kappa_B=1$, and initial conditions $R_i=B_i=1$. }
\label{largeopt2}
\end{center}
\end{figure}

\begin{figure}[tbp]
 \begin{center}
  \includegraphics[width=.99\textwidth]{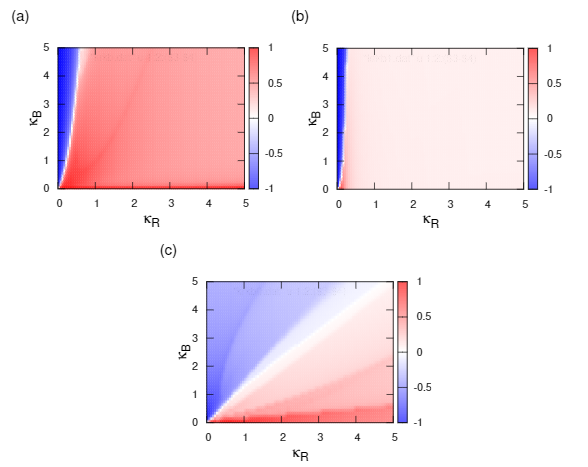}
 \caption {Heat map of differences between normalised stationary Red and Blue forces vs. lethalities $\kappa_R$ and $\kappa_B$ for the best configuration found for networks optimized for $\kappa_R=0.5, \kappa_B=1$ and (a) $\lambda=0.5$ and (b) $\lambda=0.9$. Panel (c) shows the heatmap for a random 
 Blue and Red network with the numbers of links in the manoeuvre and engagement networks consistent with the optimised networks of Red. Red colours indicate regimes of larger Red surviving forces and Blue colours regimes of larger Blue surviving forces. In all cases the system size is $N_R=N_B=50$, $L_R=L_B=100$, $\kappa_R=0.5, \kappa_B=1$, and initial conditions $R_i=B_i=1$. }
\label{largeopt2a}
\end{center}
 \end{figure}

Now we scrutinise the optimisation outcomes for a system with
numbers of nodes $N_R=N_B=50$, numbers of manoeuvre links $L_R=L_B=100$, and numbers of engagement links $L_{RB}=10$, and rates $\kappa_R=0.5, \kappa_B=1$ (thus Red is inferior to Blue in fire-power), and initial conditions $R_i=B_i=1$. The Red force optimises both its manoeuvre network and the engagement. This scenario will remain the focus
for the remainder of the paper.
Fig.\ref{largeopt2} shows results for different values
of $\lambda$ exploring the offence-defence trade-offs.
The first row has $\lambda=0.2$, for emphasis
on offence, the middle $\lambda=0.5$, and the bottom
row $\lambda=0.9$, with emphasis on defence.
The left column shows the density of surviving agents for each side
and the value of the trade-off utility function
as a function of the number of iterations of the optimisation process;
the right column shows a typical network diagram with Blue and Red coloured
manoeuvre networks at the end of the optimisation and green the
symmetric engagement network; and the panels in Fig. \ref{largeopt2a} shows
the surviving densities for networks optimized at $\kappa_R=0.5$ and $\kappa_B=1$ as a function of Red's and Blue's lethalities.

In the first column we observe the (approximate) plateauing of the 
curves showing equilibration of the optimisation, 
with the utility reaching a maximum value.
We see the intuitive behaviour that,
as $\lambda$ increases with more emphasis on self-preservation,
the equilibrium value of Blue ceases to be close to zero (panel (a)) but to 
a value of approximately 0.9 (panel (c)), 
though still less than the optimising Red force
which suffers only very small losses.

In the networks in the right column
(panels (d-f)) we show the structure at equilibrium and indicate the
resulting force strength at each node by its size. The difference in structure for different $\lambda$ is
noteworthy. 
For low and intermediate $\lambda$ 
we see a concentration of Red nodes that are multiply
engaged with Blue, highlighted by a black circle. Significantly, as closer inspection shows, these are very weak nodes.
These multiply engaged nodes may be termed `sacrificial
attackers', namely nodes attack many adversary nodes simultaneously. On the other hand these same nodes have few 
manoeuvre (red) links within their force, so that they are relatively
unsupported. In this sense, they are sacrificed -- many concurrent
attacks with no local support so that they are allowed to drain
of resource.
We also observe at small $\lambda$ in
panel (d) an absence
of Red hubs for manoeuvre, with a hint of clustering
of Red nodes 
for $\lambda=0.5$ in panel (e).
Thus Red's emphasis on Blue's destruction manifests
as little internal structure for its manoeuvre
but a significant structure in the formation of
sacrificial hubs. 
At high $\lambda$, with survival more significant in the 
utility function, a manoeuvre hub
becomes more prominent in panel (f), highlighted by a black diamond. As can be seen, this node does not attack Blue. The large
Red nodes on the ring boundary attack, but 
accumulate in different groups onto a few Blue targets, indicated
by black triangles.
Evidently, attacks are focused on a few Blue nodes which are quickly extinguished. These behaviours are a key insight from this approach to which we return later.

The panels of Fig. \ref{largeopt2a} shows the lethality dependence of the difference in the surviving densities for networks optimised at $\kappa_R$ and $\lambda=0.5, 0.9$ (note that the heatmap for $\lambda=0.2$ is similar to that for $\lambda=0.5$ and hence omitted),
and then run simulations across different values of $\kappa_R$ and $\kappa_B$. Such `battle outcome' heatmaps represent
a continuum of scenarios within the broad case
of near-peer adversaries. As a reference case, panel (i) 
corresponds to random networks
that are different instances of the same
probability distribution. It is obvious that a region of larger surviving Blue forces is separated from a region of larger surviving Red forces by the diagonal line $\kappa_R=\kappa_B$. Hence, given random structure of engagement and manoeuvre networks, the force with larger lethality `wins' and -- visualised by more intense hues of the respective colours-- victories are the more expressed the larger the distance to the diagonal line.

This scenario is dramatically altered for the optimised networks. In panel (a) we see the heatmap for a network optimised for $\lambda=0.5$ . It is now obvious that optimisation has dramatically reduced the area in parameter space that corresponds to Blue victories. Red victories typically correspond to a dark Red hue; as we have seen in Fig. \ref{largeopt2}(b,e), Red tends to win by eliminating the Blue force. Similarly, in panel (c) a heatmap for a configuration optimised for $\lambda=0.9$ is shown. Inspecting panel (c) a further reduction in the Blue area is apparent: Blue wins are now restricted to a very small area of parameter space in which Blue has much larger lethality than Red. However, as a trade-off, the Red victory becomes coloured in a less intense red. Here Red is no longer able to eliminate most Blue but wins by preserving most of its force and destroying a small proportion of Blue (see also Fig. \ref{largeopt2}(c,f)).

Analysis of optimisation results for particular parameter configurations above hints to the existence of different regimes for optimal configurations. In the next section we proceed with a more systematic analysis of the dependence of properties of optimal configurations on the trade-off parameter $\lambda$.

%Turning to the bottom row of Fig.\ref{largeopt2} showing the lethality dependence of the surviving
%densities, we take the networks optimised at the particular $\kappa_R$,
%and then run simulations across different values of $\kappa_R$.
%At low $\lambda$ with emphasis on Blue destruction,
%in panel (g), there is a clear swapping from Blue to Red
%dominance at $\kappa_R=0.5$ where the networks were optimised.
%However, once self-preservation of Blue becomes more important in
%the objective function there is a diminished sensitivity to the value of $\kappa_R$
%at which the optimisation was performed,
%as seen in panel (i).
%Comparing these curves to the same results using random graphs (in grey) for
%Red or Blue emphasises the importance of the optimisation in the outcome.

\subsection{Tradeoffs in optimal networks}
\label{tradeoff}

\begin{figure}[tbp]
\begin{center}
\includegraphics[width=.99\textwidth]{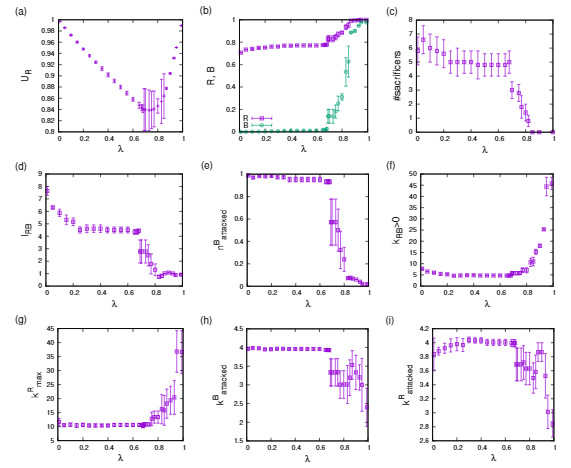}
\caption {Dependence of properties of the optimised configurations on the trade-off parameter $\lambda$ when optimising Red's utility $U_R=\lambda R+(1-\lambda) (1-B)$ (i.e. for $\lambda=0$ R is only interested in extinguishing B whereas for $\lambda=1$ force R is only interested in force preservation.) (a) Utility, (b) remaining stationary forces. (c) the number of sacrificial nodes, namely nodes
with no input in the manoeuvre network but 
engagement 
degree $k>10$, 
 (d) evolved number of links per node (attacks), 
  (e) fraction of attacked B nodes, (f) average number of attacks on attacked B-nodes, (g) evolved max. degree of R nodes in R's manoeuvre network (note that the avg. maximum degree for a random network with the same number of connections is $8.9\pm0.5$), (h) avg. degree in manoeuvre network of B nodes that are attacked by R, (i) avg. degree in manoeuvre network of R nodes that are attacking nodes of B. System of size $N_R=N_B=50$, $L_R=L_B=100$, $\kappa_R=0.5, \kappa_B=1$ (note that R is inferior to B), and initial conditions $R_i=B_i=1$. Data points represent averages over the five best optimised configurations found for the respective value of $\lambda$, if not otherwise indicated error bars are about the size of the symbols. Force R optimises both its manoeuvre
 and engagement networks.}
\label{largeopt3}
\end{center}
 \end{figure}

Now we study the dependence of the behaviours
in optimised networks against the trade-off parameter $\lambda$.
%We examine a scenario with $N_R=N_B=50$, 
%$L_R=L_B=100$, $L_{RB}=10, \kappa_B=1$ and
%$\kappa_R=0.5$. The set-up is ostensibly biased in favour of the Blue force.
Here, for every value of $\lambda$ we run the optimisation
and take the five best configurations and average properties
over these. We plot in Fig.\ref{largeopt3} the dependence of $\lambda$ of
averages of a range of quantities - the utility, numbers of evolved links, numbers of sacrificers (as discussed above),
the numbers of attacked nodes, and network degrees of various types.
Note that Red wins in every case -
as the optimised force - but the $\lambda$ dependence
of these quantities reveals different
properties for Red to achieve its
optimal performance.

Firstly, the utility function in panel (a)
shows a steady decrease as $\lambda$ increases
until approximately $\lambda=0.7$ at which
point it begins to increase again.
Thus as the two objectives - destruction of the
adversary and survival of own forces - compete
with each other the overall utility available
decreases. We see
distinct `offensive' and `defensive' phases
which swap at $\lambda=0.7$, and diminishing utility
in the intermediate regime; `pure' offence or defence objectives
($\lambda=0$ and $\lambda=1$ respectively) have higher utilities.
We emphasise that `defence' here 
always involves offence 
as fire is continuous, 
but the degree of seeking to preserve forces at the same time.
This same point of transition, $\lambda=0.7$,
manifests in all the other measures.
In the surviving forces in panel (b) we see it in
a transition from complete Blue destruction in the
offensive phase 
(it is Red that is optimised) to a coexistence
of Blue and Red for $\lambda>0.7$.

In panel (c) we compute the number of sacrificer
nodes identified as those with no
manoeuvre links with their own force but solely engaged
with a sufficiently large number of adversary nodes. We choose
$k=10$ as the threshold. We observe three phases:
a slight decrease for low $\lambda$, a plateau
and then significant drop at high $\lambda$,
again at the threshold value of $\lambda=0.7$.
In other words, sacrificial nodes are
a significant structure in the optimal networks
in the offensive phase.

The number of links per node in the evolved engagement
networks, denoted $l_{RB}$, is shown in panel (d).
Recall that the seed networks start with
$L_{RB}=10$ total number of links and 50 nodes, so $l_{RB}$, the ratio,
may change through iterations of the optimisation.
The plot in (d) 
exhibits several changes with increasing $\lambda$,
though the transition at $\lambda=0.7$
is also evident. Again, as the emphasis shifts
from offence to defence
there is a drop in the number of optimal engagement
links. This transition is sharper in the
fraction of attacked Blue nodes in panel (e):
Red goes from attacking the entire Blue force
in the offensive phase, to very few at
the critical value of $\lambda$.
In panel (f) we show
the number of attacks on attacked Blue nodes,
in other words the number of Red force agents that {\it simultaneously} attack a Blue node

This goes from a steady value
of approximately $5 \leq k_{RB} \leq 10$
to an order of magnitude larger at the critical
value of $\lambda$. In other words,
as utility shifts in emphasis to self-preservation, Red must concentrate rather
than disperse its engagement and thus chooses to engage very few Blue nodes (as seen in
(e)) but with overwhelming force. A similar
result is seen in panel (g) in the 
maximum degree of manoeuvre nodes in Red.
These nodes of high degree are the 
manoeuvre hubs
observed in the network diagrams of
Fig.\ref{largeopt2}.

To address the question whether Red attacks well supported or unsupported Blue nodes, in panels (h) and (i) we examine the
overlap between manoeuvre and engagement
in the force elements. Specifically, we measure the average manoeuvre degree of attacked Blue nodes (panel (h)) and the average degree of an attacking Red node (panel (i)).

Panel (h)
shows how Red targets the manoeuvre network of
Blue: up to the critical value of $\lambda$ almost all Blue nodes are attacked and hence the
average manoeuvre degree of 
{\it attacked} Blue corresponds to the average manoeuvre degree $k=4$, after
which it drops, though the error bars here are
quite large. In the network diagram of
Fig.\ref{largeopt2}(f) we see such
overwhelming attacks in the green links
from multiple Red nodes on to three
separate single
Blue nodes inside the ring.
The left-most one is evidently a poorly supported
Blue node, while the right-most node appears to have
a number of Blue manoeuvre links feeding into it.
Given that Blue remains an Erdos-Renyi network
(it has not been subject to optimisation)
we cannot speak here of a 
corresponding manoeuvre hub for Blue.
So, for this form of structure,
statistically in the defensive phase Red targets poorly supported Blue nodes.

In panel (i) we see the coincidence in Red of manoeuvre links
that are attacking Blue, and see a drop at the critical $\lambda$,
again with some noise. In the defensive phase Red tends to attack with nodes of low manoeuvre degree, but resources tend to be funnelled into these nodes from hub nodes in the Red manoeuvre network.
%Thus, in the defensive phase Red dedicates more nodes to a pure manoeuvre role.

 %In Fig.\ref{largeopt4} we repeat this analysis but with lower Red fire-power, $\kappa_R=0.1$. Many of these observations are repeated: the critical point remains approximately at $\lambda=0.7$ in panel (a). In panel (b) there is an even sharper transition in the number of Blue force surviving at this point when Red goes into the defensive phase. Notably in panel (g) the transition is sharper than for the higher $\kappa_R$ in the corresponding panel in Fig.\ref{largeopt3}: Red more dramatically switches away from suppressing Blue manoeuvre when it goes into the defensive state. In panel (h) we see a lower average degree of Red manoeuvre nodes that are also attacking Blue nodes (because of Red's lower lethality), with again a drop at the critical point. 
 
 We see overall a consistent pattern that matches intuition,
 on the one hand, but is non-trivial in other respects:
 the offensive phase persists beyond the trade-off point of $\lambda=0.5$,
 there is a pattern of sacrificial nodes appearing in the intermediate
 region, and manoeuvre hubs developing in the defensive regime.
 In the Appendix we show analytically how this phenomenon occurs in the model.

In the supplement we consider the case where the number of attacks
may be constrained, arguably a more realistic scenario
in warfare, and find similar behaviours to those described here
except that in the offensive phase
we find Red concentrating attacks to 
suppress the Blue force capacity to manoeuvre
resources.

\subsection{Network structure dependence on Red fire-power} 
Thus far we have fixed the fire-power, or lethality, of 
Red to a discrete set of values. We now examine the properties of
optimised networks for varying Red lethality.

For this purpose we focus on the offensive regime and fix a trade-off parameter of $\lambda=0.5$ (and $\kappa_B=1$). We then optimise networks for varying values of $\kappa_R$ and investigate their properties.

\begin{figure}[tbp]
 \begin{center}
 \includegraphics[width=.99\textwidth]{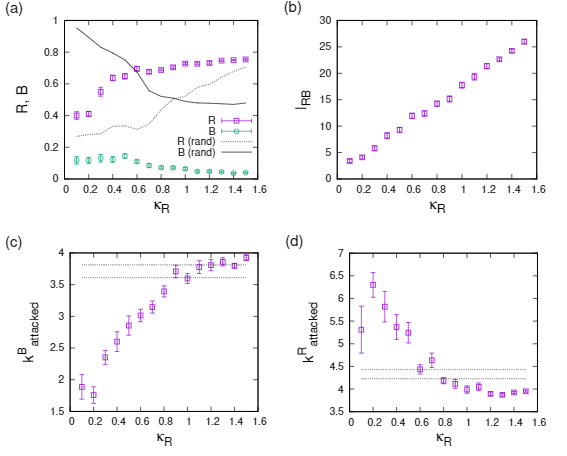}
 \caption {Outcomes for varying strength $\kappa_R$ of Red
 when manoeuvre and engagement networks are optimised for $\lambda=0.5$. (a) optimised stationary forces (solid lines give averages for a random networks that were used as seeds for the optimisation for comparison).  
 %(b) Number of attacks per node. 
 %(c) Avg. number of attacks an attacked B node receives. 
 %(d) Standard deviation of the avg. number of attacks a B node receives when attacked. 
 (b) Avg. number of B nodes under attack from R. 
 %(f) Avg. max. degree of the manoeuvre network of R. 
 (c) avg. degree in the manoeuvre network B of B nodes that are attacked by R. (d) degree in the manoeuvre network of R of R nodes that attack B.  In (c) and (d) the dotted lines indicate expectations (plus/minus one standard deviation) for random attacks calculated from the (random) networks that were used to seed the optimisation. Data points represent averages over 20 optimisation runs.}
\label{largeopt5}
\end{center}
 \end{figure}

These results are shown in the various panels of
Fig.\ref{largeopt5}. In panel (a) we 
see 
the surviving densities of the two forces with 
the result that Red is superior
for all $\kappa_R$ but with a sharp increase
between $0.3\leq \kappa_R\leq 0.5$, and
then a plateau for larger $\kappa_R$.
Thus, given the trade-off
between destruction of the adversary and self-preservation at this value of $\lambda$, there
are diminishing returns for Red increasing its lethality beyond $\kappa_R=0.5$. 
There is a corresponding flattening of
Blue's curve at large $\kappa_R$. We superimpose on
these plots the results for the original random Blue and Red graphs
which seeded the optimisation procedure with the clear result
that Red only overcomes Blue when it has the same lethality
in the absence of optimal networks.
Thus the optimisation for Red genuinely improves its performance
across all values of $\kappa_R$ compared to
the initial seed configurations.

%Panel (b), showing the number of attacks undertaken per node of Red, further exhibits this diminishment of return in increasing $\kappa_R$ where the increase tapers off around $\kappa_R=0.5$. Moreover, in panel (c) the average number of attacks received per attacked node of Blue increases then decreases beyond $\kappa_R=0.5$. The standard deviation  in (d) confirms the variability of outcomes in the low $\kappa_R$ regime. 

In (b) we show the number of attacks per node of Blue,
namely the degree of Blue nodes in the engagement
network. This evidently 
linearly increases with $\kappa_R$. A more powerful Red force will carry out more attacks than a less powerful Red force.

In
(c) we give the average degree of manoeuvre
nodes in Blue that are attacked by Red
which shows a linear increase then at $\kappa_R=1$ flattening to the average degree of the overall 
Blue manoeuvre network, $k=4$.
For low lethality Red cannot overwhelm the entire Blue force. Thus, in this regime
Red focuses its attack on selected Blue nodes which are precisely those which are least supported by Blue's manoeuvre network. For increasing lethality Red then targets most Blue nodes, thus explaining the plateau.

Finally, (d) reveals how much of Red's network is
`dual-purpose', manoeuvre and engagement, 
by computing the manoeuvre degree of Red nodes
that also attack Blue. The result is
that at $\kappa_R \approx 1$
there is again, a transition: for low fire-power
there are high-degree Red manoeuvre hubs that
also engage Blue while at high fire-power
these dual-purpose hubs have disappeared.
In light of the evidence for 
Red manoeuvre hubs in the networks of Fig.\ref{largeopt2},
this suggests that at high fire-power the hubs
become more specialised, purely manoeuvre or attack.

In (c) and (d) we also give in the dotted lines the expectation within
one standard deviation of the result using the original
random networks that seeded the optimisation. We
see that the results are constant in $\kappa_R$,
so that there is nothing distinguishing in 
the networks between manoeuvre and engagement nodes
in the overall combat outcome, due to the
uniform random nature of the Erdos-Renyi graph.
In other words, the optimisation is generating
specialised roles in the manoeuvre and engagement
structure of Red that are different at the
various values of lethality.

%%%%%%%%%%%

In summary to the analyses thus far, we conclude that
optimising with higher fire-power, when
offence and defence
are similarly weighted in utility, gives
diminishing returns. An optimised Red force structure
with lesser lethality than Blue is able to achieve
victory through dual-purposing manoeuvre and attack hubs. 
In the supplement, when the optimal force is constrained
in the number of targets it may engage,
in the offensive phase it focuses on the adversary ability to manoeuvre
its own forces. We shall discuss later the extent to which
these results resonate with ideas in Manoeuvre Warfare.

\subsection{Battle-outcome heatmaps for optimised networks}
We now take exemplar optimal networks and scan across
the parameters for the Red force, $\gamma_R, \kappa_R$,
to examine the regions where that force has
advantage over Blue.
This type of analysis would be typical as an application of
our model for trade-offs in investment
over technology that either enhances fire-power or
manoeuverability of a force but also where optimisation
of structure, both for manoeuvre and engagement,
is a consideration.

\begin{figure}[tbp]
 \begin{center}
 \includegraphics[width=.99\textwidth]{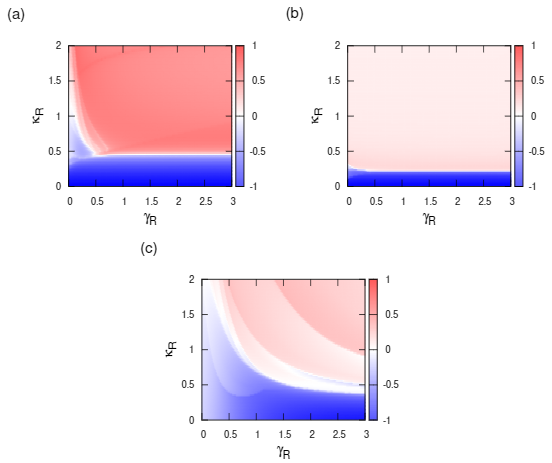}
 \caption {Plots of the victory heat-map showing the difference of 
 steady-state densities for Red against Blue across values of $\gamma_R,\kappa_R$
 for $N_B=N_R=50,L_B=L_R=100$, and networks optimised
 for $\kappa_R=0.5,\gamma_R=1,\kappa_B=1, \gamma_B=1$
 but solving the system for $\kappa_B=3$.
 Panel (a) for $\lambda=0.5$, (b) for $\lambda=0.9$, and
 (c) the result using the original random network that seeded
 the optimisation.}
\label{largeopt6}
\end{center}
 \end{figure}

We consider again forces of equal size
$N=50$ and number of manoeuvre links $L_R=L_B=100$
and optimise the Red force at 
$\kappa_R=0.5,\gamma_R=1$ against
a Blue force with $\kappa_B=1, \gamma_B=1$.
We examine two cases, with trade-off
parameter $\lambda=0.5$ for the offensive phase and $\lambda=0.9$ for the defensive phase.
We then solve the system to determine
steady-state densities of force at $\kappa_B=3$, for a 
{\it more powerful adversary than the one
for which Red has optimised}. In Fig.\ref{largeopt6}
we plot battle-outcome heatmaps,
similar to those in Figs.\ref{largeopt2} (g-i) but for Red's variable
choices, $\gamma_R,\kappa_R$.

%%%%%%%%%%%
We note the following key features in these
results in Fig.\ref{largeopt6}. Firstly, in all cases of optimised network
there is a threshold lethality required before
Red has an advantage over Blue. In panel (a), this
threshold is close to that used for 
the optimisation, $\kappa_R=0.5$, however
in panel (b), where 
self-preservation is emphasised, the critical value is even less
than that used for optimisation. We note
that these are significantly less than the
actual value of lethality for the Blue force; the optimisation is robust against variation
in adversary lethality. In both cases,
the threshold lethality is {\it insensitive} to
the manoeuvre rate $\gamma_R$.

The exception to this insensitivity is
at low $\gamma_R$ for panel (a) where $\lambda=0.5$:
here the Blue force has an advantage 
for $\gamma_R\leq 0.5$ unless the
Red force increases its lethality $\kappa_R$.
Contrastingly, in panel (b) the Red force maintains advantage
for $\kappa_R\geq 0.25$, 
right down to very small values of $\gamma_R$.
Recalling the key result from our previous analyses,
when offence is more valued in utility
for $\lambda\leq 0.7$, the optimised Red networks
strongly figure sacrificer nodes as a key 
mechanism in the engagement structure. As we saw in previous analyses, for
$\lambda>0.7$, these disappear and are replaced
%AK: tried changing to get both mechanisms in
by Red manoeuvre hubs while all the remaining Red
nodes engage with a small number of Blue. Thus in Fig.\ref{largeopt6}(a),
the optimised network will have more emphasis
on attack rather than manoeuvre, with this
flipping for panel (b). It is clear then that
for (a), at low $\gamma_R$ the optimised
manoeuvre network is insufficient to
guarantee victory - additional lethality is required.
Contrastingly, in (b) the {\it manoeuvre} network
gains prominence in the optimisation through the manifestation of
hubs while the rest of Red are
engaged in concentrating attacks on Blue.
Thus {\it fewer} Red require to manoeuvre (they purely attack) {\it except through an efficient structure}, a hub.
This leads to low $\gamma_R$ 
values being
sufficient for effectiveness against the Blue force.

Comparing Figs.\ref{largeopt6}(a,b) we also see that
in the region where Red does have advantage
it scores {\it better} when {\it equally} weighting 
offence and defence (panel (a)) compared to
emphasis on defence (in (b)). This seems in
apparent contradiction to Fig.\ref{largeopt3}(b) where
at high $\lambda$ the number of surviving Red forces is larger
than at lower $\lambda$. However, at lower $\lambda$
{\it no} Blue forces survive, unlike the case for larger $\lambda$.
Thus the {\it difference} between steady-state force numbers 
at $\lambda=0.9$ may indeed be less than for $\lambda=0.5$.

In Fig.\ref{largeopt6}(c) we show the result
using the original random networks that seeded
the optimisation; both Red and Blue have Erdos-Renyi random graphs. 
The shape of the cross-over region is
understandable: Red has advantage when it is high in
manoeuvre rate $\gamma_R$ and fire-power $\kappa_R$.
Contrasting this with the
optimised results, we see the sensitivity in
both of Red's parameters with gradients even beyond the threshold
for Red victory. Thus the key role of optimisation is in flattening
the variability in the various parameters.

To summarise these results, we see that
with the capacity to invest in
optimised networks, and fire-power or speed of
manoeuvre through its structures,
the Red force gains advantage through
allocation of resources {\it up to
a minimum level in fire-power}. Beyond that,
it is wasted resource. Similarly,
as long as it has optimised structure for
both manoeuvre and engagement, the force only needs
a minimal threshold in manoeuvre rate to gain advantage.

\section{Conclusions and Discussion}
We have generalised the $(N,M)$ Lanchester model 
to include manoeuvre warfare through the incorporation of networks. 
The model allows the embedding of warfighting
heuristics for how a force dynamically redistributes
resources through battle according to 
some local weighting across connections; here
we used the local
ratio of a force element to its adversaries.
We optimised the structure of networked
forces for both manoeuvre
and engagement against a randomly structured adversary,
as functions of an offence-defence trade-off and fire-power.
We found consistently, that a force
with optimised networks
could defeat an equally sized opponent for lesser fire-power than that of
the adversary.
Depending on how it valued the offence-defence trade-off,
different structures emerged for the optimised
force.

We now reflect on how these structures figure in the 
literature of Manoeuvre Warfare.

Arguably, our model only narrowly reflects the manoeuvre dimension
of warfare, purely in terms of the movement of resources
around a networked force rather than
the decision process. For something that
is defined as {\it one} of the warfighting functions, military writings on Manoeuvre Theory tend to wrap
numerous functions together, as
well as the psychological dimensions of surprise and
shock.
Thus \cite{Lind1985} draws upon the Manoeuvre Theory
of John Boyd 
which treats Manoeuvre as a competition between adversarial Observe-Orient-Decide-Act (OODA) loops, itself a
a model for C2. 
Thus our the manoeuvre rate constant $\gamma$ 
is a proxy for the speed of each agent's OODA loop 
in the absence of an explicit structure
for distributed decision making (also enabled by networks) that may facilitate
or undermine the capacity to manoeuvre.
In this interpretation, $\gamma$ measures not just
technology (speed and agility of vehicles or computers)
but cognitive and psychological capacity.

Many authors see Sun Tzu \cite{Smith&LeBrun1994} as
the father of Manoeuvre Theory in his maxim ``To subdue the
enemy without fighting'' (III.6). Clearly in our model,
fighting is continuous through the process represented
by the differential equations.
However, implementations of Manoeuvre Theory 
into US Army doctrine, such as 
the 1970-1990s US Defense Reformers who emphasised
(some argue, erroneously\cite{Lauer1991}) 
``smaller, more mobile forces'' against larger
adversaries, {\it preserving one's own
overall resource  during the fight} is the real point of
Sun Tzu. This our model captures and leads to recognisable structures. 
To this end we may focus on
Clausewitz's {\it Principles of War} \cite{Clausewitz1942},
the tactical precursor to his most well-known work
{\it On War} \cite{Clausewitz1976}:
``We must select for our attack one point of the enemy's position and attack it with great superiority, leaving the rest of his army in uncertainty but keeping it occupied. This is the only way that we can use an equal or smaller force to fight with advantage and thus with a chance of success. The weaker we are, the fewer troops we should use to keep the enemy occupied at unimportant points, in order to be as strong as possible at the decisive point.'' (Sect II.2.1).
There follows a set of diagrams outlining how
reserves may be used in certain tactical configurations,
with forces circling an enemy and reserves in interior lines.

Though often contrasted with Sun Tzu as an attritionist, 
one can recognise here that Clausewitz
espouses mobility of weaker forces
to achieve {\it localised concentration} on the enemy.
While valuing decisive battle over Sun Tzu's elusive force, 
Clausewitz emphasise manoeuvre in getting to that fight.

Analysing these ideas in reverse we see first in
the manoeuvre hubs arising in the model for $\lambda\approx 1$ a pattern
of how reserves should flow into battle, abstracted from
spatial location, identifiable with
Clausewitz's ``formation in depth''.
The sacrificial nodes that emerged in the intermediate
region of the offence-defence trade-off $\lambda$ may be interpreted as {\it feints}
-- ``unimportant points" --
insofar as they distract the adversary fire with
`mass scale expendable' resources.
These may be `dummy' units, as
we alluded in the introduction,
which are replete in military history.
Successful feints are based on incomplete 
information - 
the enemy lacks knowledge of the true nature of
these targets. In
our model the
lack of information is reflected in the fact that 
engagement targets (of Red, from the Blue force) are assigned non-adaptively;
only the manoeuvre
of resource given the state of the fight is adaptive. If we were to
use more sophisticated adaptive
engagement mechanisms, such as
examined by \cite{MacKay2009} we can
anticipate that this will change. The ignorance implied
for such an interpretation of the sacrificial
node behaviour is consistent with
the absence of a representation of ISR
at this stage of the model development.

The point of concentration of fire on the enemy is
often described as the {\it schwerpunkt}.
In our analyses, Blue is not optimised and thus
is not adaptive to form structural or spatial hubs, thus offering Red a functional {\it schwerpunkt} for its attacks.  
The force with optimised networks thus attacks the
ability of its adversary to manoeuvre resources,
particularly when the number of attacks
is constrained (as in the supplement), 
through what may be interpreted as
{\it suppressive fires}. 
Finally, we may recognise
in the collocation of manoeuvre and engagement hubs 
for the optimised force qualities seen in
the coincidence of armour, speed, fire-power
and mobile communications of the {\it panzer} divisions
from 1940-41 of World War II. Indeed, this
form of warfare was the realisation of many of the visions
of manoeuvrist theorists, both British and German
included, in the aftermath of the Great War.

These are qualitative observations
about the patterns emerging in this mathematical
approach. As alluded from the outset,
the value of this approach is that it may be
readily generalised further to include
alternate heuristics for resource and target reallocation, and 
further warfighting functions.
In particular, current work is developing this
model to incorporate C2 through 
models of synchronising dynamical processes on network,
effectively modelling the interactions of multiple
Boyd OODA loops in the elements of each force. 
Logistics, built on yet other networks, are also straightforwardly
incorporated into the model as source terms to represent how
forces enter into and sustained through the battle.
Also, using various differential equation based models of swarming,
spatially embedding the model is also straightforward.
It would be at this stage that
appropriate validation studies should be undertaken
to compare with data from a range of manoeuvrist
centric battles in history. Nevertheless, we would agree
with MacKay that the value of such models
is less to be predictive (for which
the quality of data is often lacking), but more to explore
more sophisticated warfighting concepts
such as the role of, and balance of investments,
across competing technologies, doctrine  
and human training to realise the various warfighting functions.

\section*{Acknowledgements}

We are grateful for discussions with Ryan Ahern, Sharon Boswell and Brandon Pincombe. This work was conducted under the auspices of DST's Modelling Complex Warfighting Initiative.
MB acknowledges partial support from IES-R2-192206 and the Alan Turing pilot project ID 519465.
The simulations reported here were run on the supercomputer IRIDIS5 at the University of Southampton.

%\section*{Appendices}

%Any appendices should be placed after the list of references, beginning with the command \verb"\appendix" followed by the command \verb"\section" for each appendix title, e.g.

\section*{Appendix: An analytical explanation of `sacrificer' nodes}

%\noindent\textbf{Appendix. An analytical explanation of `sacrificer' nodes}\medskip
\begin{figure}[tbp]
 \begin{center}
 \includegraphics[width=.5\textwidth]{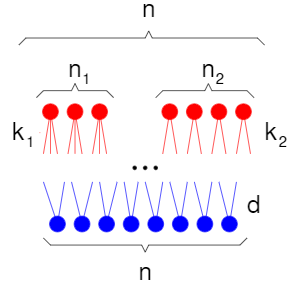}
 \caption {Illustration of a red force of $n$ nodes comprised of two subgroups of size $n_1$ and $n_2$ who, respectively, each carry out $k_1$ and $k_2$ attacks on a blue force of size $n$, assuming that each blue node is attacked by $d$ red nodes.}
\label{sacrifs}
\end{center}
 \end{figure}
 
 In subsection 5.3, we saw that when optimising the combined engagement-manoeuvre networks a regime for low $\lambda$ exists in which optimal configurations show the coexistence of `normal' Red nodes that attack one Blue node each and a certain number of `sacrificial' Red nodes which typically attack many Blue nodes, but don't receive support via the Red manoeuvre network. Here, we will show analytically that such nodes are indeed typical in our present setup, even without the presence of manoeuvre support. Consider the set-up depicted in Fig. \ref{sacrifs}, in which the Red force composed of $n$ nodes consists of two groups $n_1+n_2=n$ and fights a Blue adversary force made up of $n$ nodes. Suppose, Red nodes in group one each attack $k_1$ Blue nodes and Red nodes in group two each attack $k_2$ nodes. Below, we will argue about parameters $k_1$, $k_2$, $n_1$, and $n_2$ that optimise outcomes for Red. For simplicity, we will further assume that Red nodes target Blue nodes at random. As there is a total of $L=L_1+L_2=n_1 k_1+n_2 k_2$ attacks from the Red force, on average a Blue node will be targeted by $d=L/n$  Red nodes. Of these attacks, $L_1/L$ are attacks from group one, and $L_2/L$ attacks from group two. Overall, treating all Red nodes within a group and Blue nodes each as an average representative node, we thus arrive at the following system of equations governing the `mean-field' Lanchester dynamics
\begin{align}
 \label{ELan1} \dot{R_1} &= -\kappa_B k_1 B \frac{n}{L}\\
 \label{ELan2} \dot{R_2} &= -\kappa_B k_2 B \frac{n}{L}\\
 \label{ELan3} \dot{B}  &= -\kappa_R \frac{L_1}{n} \frac{1}{k_1} R_1-\kappa_R \frac{L_2}{n} \frac{1}{k_2} R_2,
\end{align}
where $R_1$ ($R_2$) is the force at an average group one (group two) Red node, and $B$ the force of an average Blue node. Respectively multiplying Eqs. (\ref{ELan1}) and (\ref{ELan2}) by $-\kappa_R L (L_1/n^2)(R_1/k_1^2)$ and $-\kappa_R L (L_2/n^2)(R_2/k_2^2)$ and Eq. (\ref{ELan3}) by $\kappa_B B$ gives
\begin{equation}
    \frac{d}{dt} \left( -\kappa_R \frac{L}{n^2} \left( \frac{n_1}{k_1}R_1^2 + \frac{n_2}{k_2} R_2^2 \right) + \kappa_B B^2\right) = 0,
\end{equation}
 and we have thus found a typical invariant for the Lanchester dynamics. Suppose that, as also in our computational experiments, initial allocations of Blue and Red forces are equal. Then, in order to maximize battle outcomes for Red, the Red commander has to maximize the function
 \begin{align}
    f(k_1,k_2,n_1)= \frac{L}{n^2} \left( \frac{n_1}{k_1} + \frac{n_2}{k_2}\right),
 \end{align}
 where $n_2=n-n_1 \geq 0$. Straightforward analysis then shows that $f$ is maximized for $n_1=n/2$ and $k_1=1$ and $k_2=n$ (or $k_1=n$, $k_2=1$). Consequently, optimal results for Red are achieved if Red splits its force into one group of nodes that each attack all Blue nodes and one group of nodes who each attack exactly one Blue node, namely if Red diverts a substantial amount of its forces as sacrificial nodes.
 
 The advantage Red can gain by introducing sacrificial nodes is seen when we evaluate comparable invariants of the Lanchester dynamics. We find the following victory conditions for Red. Presuming that each Red node initially has equal force $R_1(t=0)=R_2(t=0)=R_0$ and $B(t=0)=B_0$,  at the optimal engagement configuration we have
 \begin{equation}
     \kappa_R R_0^2 \left( \frac{1}{2} + \frac{n}{4} + \frac{1}{4n} \right) > \kappa_B B_0^2,
 \end{equation}
 whereas for a non-optimised engagement in which $k=k_1=k_2$ one has
 \begin{equation}
     \kappa_R R_0^2 > \kappa_B B_0^2.
 \end{equation}
 Thus, in the limit of very large forces $n\gg 1$ a Red commander who optimises its engagement structure only needs to bring a fraction $2/\sqrt{n}$ of the force to gain victory compared to the non-optimised configuration.

%\noindent\textbf{Appendix B. This is the title of the second appendix}\medskip

%\noindent Subsections, equations, figures, tables, etc.\ within appendices will then be automatically numbered as appropriate. Some theorem-like environments may need to have their counters reset manually (e.g.\ if they are not numbered within sections in the main text). You can achieve this by using \verb"\numberwithin{remark}{section}" (for example) just after the \verb"\appendix" command.

%Please note that if the \verb"endfloat" package is used on a document containing appendices, the \verb"\processdelayedfloats" command must be included immediately before the \verb"\appendix" command in order to ensure that the floats in the main body of the text are numbered as such.

%\processdelayedfloats %%% See above for an explanation of why this command might be needed.

%\appendix

\end{document}